\documentclass[11pt]{amsart}
\usepackage{amssymb,amsmath,amsgen,amsxtra,amsfonts,a4wide} 
\usepackage{dsfont,times}

\begin{document}
%
%
%
\theoremstyle{definition}
\newtheorem{Definition}{Definition}[section]
\newtheorem{Convention}{Definition}[section]
\newtheorem{Construction}{Construction}[section]
\newtheorem{Example}[Definition]{Example}
\newtheorem{Examples}[Definition]{Examples}
\newtheorem{Remark}[Definition]{Remark}
\newtheorem{Remarks}[Definition]{Remarks}
\newtheorem{Caution}[Definition]{Caution}
\newtheorem{Conjecture}[Definition]{Conjecture}
\newtheorem{Question}[Definition]{Question}
\newtheorem{Questions}[Definition]{Questions}
\theoremstyle{plain}
\newtheorem{Theorem}[Definition]{Theorem}
\newtheorem*{Theoremx}{Theorem}
\newtheorem{Proposition}[Definition]{Proposition}
\newtheorem{Lemma}[Definition]{Lemma}
\newtheorem{Corollary}[Definition]{Corollary}
\newtheorem{Fact}[Definition]{Fact}
\newtheorem{Facts}[Definition]{Facts}
\newtheoremstyle{voiditstyle}{3pt}{3pt}{\itshape}{\parindent}%
{\bfseries}{.}{ }{\thmnote{#3}}%
\theoremstyle{voiditstyle}
\newtheorem*{VoidItalic}{}
\newtheoremstyle{voidromstyle}{3pt}{3pt}{\rm}{\parindent}%
{\bfseries}{.}{ }{\thmnote{#3}}%
\theoremstyle{voidromstyle}
\newtheorem*{VoidRoman}{}

%
\newcommand{\prf}{\par\noindent{\sc Proof.}\quad}
\newcommand{\blowup}{\rule[-3mm]{0mm}{0mm}}
\newcommand{\cal}{\mathcal}
\newcommand{\Aff}{{\mathds{A}}}
\newcommand{\BB}{{\mathds{B}}}
\newcommand{\CC}{{\mathds{C}}}
\newcommand{\FF}{{\mathds{F}}}
\newcommand{\GG}{{\mathds{G}}}
\newcommand{\HH}{{\mathds{H}}}
\newcommand{\NN}{{\mathds{N}}}
\newcommand{\ZZ}{{\mathds{Z}}}
\newcommand{\PP}{{\mathds{P}}}
\newcommand{\QQ}{{\mathds{Q}}}
\newcommand{\RR}{{\mathds{R}}}
\newcommand{\Liea}{{\mathfrak a}}
\newcommand{\Lieb}{{\mathfrak b}}
\newcommand{\Lieg}{{\mathfrak g}}
\newcommand{\Liem}{{\mathfrak m}}
\newcommand{\ideala}{{\mathfrak a}}
\newcommand{\idealb}{{\mathfrak b}}
\newcommand{\idealg}{{\mathfrak g}}
\newcommand{\idealm}{{\mathfrak m}}
\newcommand{\idealp}{{\mathfrak p}}
\newcommand{\idealq}{{\mathfrak q}}
\newcommand{\idealI}{{\cal I}}
\newcommand{\lin}{\sim}
\newcommand{\num}{\equiv}
\newcommand{\dual}{\ast}
\newcommand{\iso}{\cong}
\newcommand{\homeo}{\approx}
\newcommand{\mm}{{\mathfrak m}}
\newcommand{\pp}{{\mathfrak p}}
\newcommand{\qq}{{\mathfrak q}}
\newcommand{\rr}{{\mathfrak r}}
\newcommand{\pP}{{\mathfrak P}}
\newcommand{\qQ}{{\mathfrak Q}}
\newcommand{\rR}{{\mathfrak R}}
%
%
\newcommand{\dq}{{``}}
\newcommand{\OO}{{\cal O}}
\newcommand{\into}{{\hookrightarrow}}
\newcommand{\onto}{{\twoheadrightarrow}}
\newcommand{\Spec}{{\rm Spec}\:}
\newcommand{\Proj}{{\rm Proj}\:}
\newcommand{\Pic}{{\rm Pic }}
\newcommand{\Br}{{\rm Br}}
\newcommand{\NS}{{\rm NS}}
\newcommand{\chit}{\chi_{\rm top}}
\newcommand{\KanDiv}{{\cal K}}
\newcommand{\perdef}{{\stackrel{{\rm def}}{=}}}
\newcommand{\Cycl}[1]{{\ZZ/{#1}\ZZ}}
\newcommand{\Sym}{{\mathfrak S}}
\newcommand{\Xcan}{X_{{\rm can}}}
\newcommand{\ab}{{\rm ab}}
\newcommand{\Aut}{{\rm Aut}}
\newcommand{\Hom}{{\rm Hom}}
\newcommand{\Supp}{{\rm Supp}}
\newcommand{\ord}{{\rm ord}}
\newcommand{\divisor}{{\rm div}}
\newcommand{\Alb}{{\rm Alb}}
\newcommand{\Jac}{{\rm Jac}}
\newcommand{\piet}{{\pi_1^{\rm \acute{e}t}}}
\newcommand{\Het}[1]{{H_{\rm \acute{e}t}^{{#1}}}}
\newcommand{\Hcris}[1]{{H_{\rm cris}^{{#1}}}}
\newcommand{\HdR}[1]{{H_{\rm dR}^{{#1}}}}
\newcommand{\hdR}[1]{{h_{\rm dR}^{{#1}}}}
\newcommand{\defin}[1]{{\bf #1}}
\title[Algebraic Surfaces with Small $c_1^2$]{Algebraic Surfaces of General Type with Small $\mathbf{c}_\mathbf{1}^\mathbf{2}$\\ in Positive Characteristic}
\author{Christian Liedtke}
\address{Mathematisches Institut, Heinrich-Heine-Universit\"at, 40225
  D\"usseldorf, Germany}
\email{liedtke@math.uni-duesseldorf.de}
\date{October 26, 2007}

\begin{abstract}
We establish Noether's inequality for surfaces of general type
in positive characteristic.
Then we extend Enriques' and Horikawa's classification of surfaces
on the Noether line, the so-called Horikawa surfaces.
We construct examples for all possible
numerical invariants and 
in arbitrary characteristic, where we need
foliations and deformation techniques 
to handle characteristic $2$.
Finally, we show that Horikawa surfaces lift
to characteristic zero.
\end{abstract}
\setcounter{tocdepth}{1}
\maketitle
\tableofcontents
\section*{Introduction}

The genus $g$ and the degree of a canonical divisor $K_C$ 
of a smooth projective curve are related by the
well-known formula $\deg K_C=2g-2$.

Already in 1875,
Max~Noether \cite{no} has given the following generalisation
to surfaces:
Given a minimal surface of general type $X$ 
over the complex numbers with 
self-intersection $K_X^2$ of a canonical divisor
(playing the role of the degree of a canonical divisor)
and with geometric genus $p_g$
then
$$
K_X^2\,\geq\,2p_g\,-\,4\,.
$$
It is natural to classify surfaces for which equality holds,
i.e., surfaces on the so-called {\em Noether line}.
Over the complex numbers, this has been sketched in
Enriques' book \cite[Capitolo VIII.11]{enr} 
and a detailed analysis has
been carried out by Horikawa
\cite{horikawa}.
The result is that surfaces on the Noether line, also called
{\em Horikawa surfaces}, are double covers of
rational surfaces via their canonical map.
Hence these surfaces may be thought of as a two-dimensional
generalisation of hyperelliptic curves.

Another point of view comes from \v{S}afarevi\v{c}'s book
\cite[Chapter 6.3]{saf}:
The $3$-canonical map of a complex surface of general type is
birational as soon as $K_X^2>3$.
Hence, by Noether's inequality, surfaces with $p_g>3$ 
have a birational $3$-canonical map.
However, Horikawa surfaces with $p_g=3$ provide
examples of surfaces where the $3$-canonical map is not
birational.
Also this may be thought of as a generalisation of
hyperelliptic curves.
\medskip

In this article, we extend these results to surfaces
of general type
over algebraically closed fields of arbitrary characteristic.
We first show that Noether's inequality still holds.
Our contribution here lies in skipping through the literature
to find a characteristic-free proof.

\begin{VoidItalic}[Theorem \ref{noether}]
 Let $X$ be a minimal surface of general type. 
 Then
 $$K^2_X\,\geq\,2p_g\,-\,4.$$
 If the canonical system is composed with a pencil
 and $p_g\geq3$, then
 $K_X^2\,\geq\,2p_g\,-\,2$ holds true.
\end{VoidItalic}

\begin{VoidItalic}[Definition \ref{horikawadefinition}]
  A {\em Horikawa surface} is a minimal surface $X$ of general
  type for which the equality
  $K_X^2=2p_g\,-\,4$ holds.
\end{VoidItalic}
\medskip

In order to extend the classification of Enriques and Horikawa
of these surfaces we need Clifford's theorem on special linear 
systems for singular curves.
With this result we can avoid Bertini's theorem 
in the classical argumentation and obtain literally the
same result, now valid in all characteristics.
As usual, there is an extra twist in characteristic $2$.

\begin{VoidItalic}[Theorem \ref{classification} and Proposition \ref{restrictions}]
  Let $X$ be a Horikawa surface and $S:=\phi_1(X)$ the image of the
  canonical map, which is a possibly singular surface in $\PP^{p_g-1}$.
  
  Then $\phi_1$ is a generically finite morphism of degree $2$ and we have the
  following cases:
  \begin{enumerate}
  \item If $S$ is a smooth surface then we have the
    following possibilities
    \begin{itemize}
    \item[-] $S\iso\PP^2$ and $p_g=3$.
    \item[-] $S\iso\PP^2$ and $p_g=6$.
    \item[-] $S\iso\FF_d$ and $p_g\geq\max\{d+4, 2d-2\}$
      and $p_g-d$ is even.
    \end{itemize}
  \item If $S$ is not smooth then it is the cone over the rational normal
    curve of degree $d:=p_g-2$.
    The minimal desingularisation of $S$ is the Hirzebruch
    surface $\FF_{d}$ and
    $4\leq p_g\leq 6$.
  \end{enumerate}
\end{VoidItalic}

From this description we deduce that Horikawa surfaces
are algebraically simply connected and that their Picard
schemes are reduced, cf. Proposition \ref{simplyconnected}.
Another byproduct is 
Proposition \ref{non3birational}, which tells us
that the $3$-canonical map of a Horikawa surface with
$p_g=3$ is not birational onto its image.
\medskip

Conversely, out of this data one can always construct a
Horikawa surface.
In characteristic $p\neq2$ this can be done along the lines
of Horikawa's article \cite[Section 1]{horikawa}, cf.
Section \ref{inseparablesection}.
\medskip

Hence we are interested in Horikawa
surfaces in characteristic $2$.
To obtain such surfaces, we use
quotients of minimal rational surfaces by $p$-closed foliations.
The main technical difficulty is that these vector fields
necessarily have isolated singularities so that we need
to control the singularities of the quotients.

The canonical map of such a surface is
a purely inseparable morphism onto a rational surface.
In particular, these surfaces are inseparably unirational.
This is in contrast to curves, whose canonical
maps are always separable.

\begin{VoidItalic}[Theorem ({\rm Horikawa \cite[Section 1]{horikawa}, %
  Section \ref{inseparablesection}})]
  All possible cases of the previous theorem do exist in
  arbitrary characteristic.
  In characteristic $2$, we may even assume the canonical
  map to be inseparable.
\end{VoidItalic}

To get Horikawa surfaces in
characteristic $2$ with separable canonical
map, we use a deformation argument.
Morally speaking, it says that surfaces with inseparable canonical
map should be at the boundary of the moduli space.
In particular, all possible numerical invariants for a Horikawa
surface in characteristic $2$ do occur with surfaces that have a
separable canonical map.

\begin{VoidItalic}[Theorem \ref{deformation}]
  In characteristic $2$, every Horikawa surface with inseparable
  canonical map can be (birationally) deformed into a Horikawa
  surface with a separable canonical map, while fixing 
  $p_g$ and the canonical image.
\end{VoidItalic}
\medskip

Since the classification looks the same in every characteristic,  
it is natural to ask whether Horikawa
surfaces over an algebraically closed field of positive characteristic $k$
lift over the Witt ring $W(k)$.
I.e., we look for a scheme $\cal X$, flat over $\Spec W(k)$ and
with special fibre the given surface over $k$.

\begin{VoidItalic}[Theorem \ref{lifting}]
  The canonical model of a Horikawa surface lifts over
  $W(k)$.
  For every Horikawa surface there exists an
  algebraic space, flat over a possibly ramified 
  extension of $W(k)$, that achieves the lifting.
\end{VoidItalic}

As a Horikawa surface is a double cover of a rational surface,
the idea of proving Theorem \ref{lifting} is first to lift the rational
surface and then the line bundle associated with this double cover, which
then defines a lifting of the whole double cover and hence a lifting
of the Horikawa surface in question.

\begin{VoidRoman}[Acknowledgements]
  I thank Stefan~Schr\"oer for many discussions and help.
  Also, I thank Matthias Sch\"utt and the referee for pointing out 
  a couple of inaccuracies. 
\end{VoidRoman}

\section{Singular hyperelliptic curves}

In order to prove Theorem \ref{firsttheorem} below, we recall
some facts about singular hyperelliptic curves.
In this section, curves will always be assumed to be
reduced and irreducible as well
as proper over an algebraically closed field $k$.
We denote by $p_a(C)\,:=\,1-\chi(\OO_C)$ the arithmetic genus
of $C$.

\begin{Definition}
  \label{hyperdef}
  A reduced and irreducible curve $C$
  is called {\em hyperelliptic} if $p_a(C)\geq2$ and if
  there exists a morphism of degree $2$ from $C$ onto $\PP^1$.   
\end{Definition}

It follows that $C$ is automatically Gorenstein, say with invertible
dualising sheaf $\omega_C$.
Clearly, a smooth curve is hyperelliptic in the sense of Definition \ref{hyperdef}
if and only if it is hyperelliptic in the classical sense.

An immediate consequence that will be used later on is the following
result.

\begin{Lemma}
  \label{nonbirational}
  Let $\cal L$ be an invertible sheaf on a hyperelliptic curve
  such that
  ${\cal L}^{\otimes 2}\iso\omega_C$.
  Then 
  $|{\cal L}|$ does not define a birational map.
\end{Lemma}

\prf
Let $\phi:C\to\PP^1$ be a morphism of degree $2$.
The fibres of $\phi$ provide us with an infinite number of
smooth points $x,y$ (possibly $x=y$) such that 
$h^0(\OO_C(x+y))\geq2$.

First, suppose that $\phi$ is generically \'etale.
Then we may assume $x\neq y$.
The long exact sequence in cohomology,
Serre duality and the assumption on $\cal L$ yield
$$
\begin{array}{cccccccc}
0&\to&H^0({\cal L}\otimes\OO_C(-x-y))&\to&H^0({\cal L})&\to&
{\cal L}/(\mm_x\oplus\mm_y)\cdot{\cal L}&\stackrel{\delta}{\to}\\
&&H^0({\cal L}\otimes\OO_C(x+y))^\vee&\to&
H^0({\cal L})^\vee&\to&0\,.
\end{array}
$$
By the choice of $x,y$ and
\cite[Lemma IV.5.5]{hart}, we have
$$
h^0({\cal L}\otimes\OO_C(x+y)) \,\geq\,
h^0({\cal L})\,+\,h^0(\OO_C(x+y))\,-\,1 \,\geq\,
h^0({\cal L})\,+\,1\,,
$$
which implies that the boundary map $\delta$
is non-trivial.
In particular, $|{\cal L}|$ fails to separate
an infinite number of distinct points
$x$ and $y$ and so the associated map 
cannot be birational.

If $\phi$ is not generically \'etale, then its fibres provide
us with an infinite number of points where $|{\cal L}|$ fails
to be an embedding (similar long exact sequence as before).
Again, $|{\cal L}|$ cannot be birational.
\qed\medskip

It is interesting that Clifford's theorem on special line
bundles remains true in the singular case.

\begin{Theorem}[Clifford's theorem]
  \label{clifford}
  Let $\cal L$ be an invertible sheaf on a reduced Gorenstein curve $C$.
  If both, $h^0(C,{\cal L})$ and $h^1(C,{\cal L})$ are non-zero then
  $$
  h^0(C,\,{\cal L})\,\leq\,\frac{1}{2}\deg{\cal L}\,+\,1\,.
  $$
  Moreover, if equality holds then ${\cal L}\iso\OO_C$ or
  ${\cal L}\iso\omega_C$ or 
  $C$ is hyperelliptic.
\end{Theorem}

\prf
\cite[Theorem A]{eks}, where even a version for torsion free
sheaves is given.
\qed\medskip

Although we do not need this result in the sequel we note that
being hyperelliptic can be rephrased in
terms of the canonical map $|\omega_C|$, just as in the
smooth case.
This is proved as in the classical case or along the lines of the
proof of Lemma \ref{nonbirational}.

\begin{Proposition}
  A Gorenstein curve $C$ with $p_a(C)\,\geq\,2$
  is hyperelliptic if and only if 
  $|\omega_C|$ is not birational.\qed
\end{Proposition}

\section{Noether's inequality and Horikawa surfaces}

The results of this section are well-known over the complex numbers.
In order to extend them to positive characteristic we run 
through the classical arguments and have to find new ones whenever
Bertini's theorem or vanishing results are used.

We start with Noether's inequality \cite[Abschnitt 11]{no}. 
Probably it is known to the experts that it holds
in arbitrary characteristic.

\begin{Theorem}[Noether's inequality]
 \label{noether}
 Let $X$ be a minimal surface of general type.
 Then
 $$K_X^2\,\geq\,2p_g\,-\,4.$$
 If the canonical system is composed with a pencil
 and $p_g\geq3$, then
 $K_X^2\,\geq\,2p_g\,-\,2$
 holds true.
\end{Theorem}

\prf
Since $K_X^2>0$ we may assume that $p_g\geq3$.
Thus, the canonical system is non-empty and either is composed
with a pencil or has a $2$-dimensional image.

If the canonical linear system is composed with a pencil
and $p_g\geq3$, 
we argue as in the proof of \cite[Theorem VII.3.1]{bhpv} 
to get  the inequality $K_X^2\geq2p_g-2$ in this case.

If the canonical map has a $2$-dimensional image,
\cite[Proposition 0.1.3 (iii)]{ek2} yields the desired inequality.
\qed\medskip

In view of this inequality it is natural to classify those
surfaces where $K_X^2$ attains the minimal value possible
given $p_g$.
Over the complex numbers, this classification
has been sketched in Enriques' book \cite[Capitolo VIII.11]{enr}
and carried out in detail by Horikawa \cite{horikawa}.

\begin{Definition}
  \label{horikawadefinition}
  A {\em Horikawa surface} is a minimal surface $X$ of general
  type for which the equality
  $K_X^2\,=\,2p_g(X)\,-\,4$ holds.
\end{Definition}

Sometimes these surfaces are referred to as 
{\em even Horikawa surfaces} as $K_X^2$ is always an even number.
Since Horikawa also classified surfaces for which $K_X^2=2p_g-3$ holds
and $K_X^2$ for such a surface is an odd number these latter surfaces
are sometimes called {\em odd Horikawa surfaces}.
However, we will only deal with even Horikawa surfaces in this article,
and so we will simply refer to them as Horikawa surfaces.
\medskip

We now establish the structure result about the canonical map
of a Horikawa surface in positive characteristic.
Although our proof is essentially the same as the original 
one by Enriques and Horikawa, we have to be a little bit 
careful applying Bertini's theorem to the canonical linear system.
In fact, the Horikawa surfaces in characteristic $2$ that
we will construct in Section \ref{inseparablesection} have the
property that a generic canonical divisor is a 
singular rational curve.

\begin{Theorem}[Enriques, Horikawa]
  \label{firsttheorem}
  Let $X$ be a Horikawa surface.
  Then $p_g\geq3$ and the canonical map $\phi_1$ 
  is basepoint-free and without fixed part.
  More precisely, $\phi_1$ is a generically finite morphism
  of degree $2$ onto a possibly singular
  surface of degree $p_g-2$ inside $\PP^{p_g-1}$.
\end{Theorem}

\prf
Since Noether's inequality is an equality and since
$K_X^2>0$ for a minimal surface of general type 
we have $p_g\geq3$.
In particular, the canonical system is not empty.
By Theorem \ref{noether}, the canonical system is not
composed with a pencil and hence the image of the
canonical map $\phi_1$ is a surface.

It follows from \cite[Proposition 0.1.2 (iii)]{ek2} that $\phi_1$
is basepoint-free and either birational or of degree $2$ onto a
(possibly singular) ruled surface.

Suppose that $\phi_1$ is birational.
Let $D$ be a canonical divisor, which we can assume
to be irreducible 
by Bertini's theorem \cite[Th\'eor\`eme I.6.10]{jou}.
Being birational, $\phi_1$ is generically unramified,
which implies that $D$ is reduced over an open
and dense subset (loc. cit.).
We may thus assume that $D$ is
a reduced and irreducible curve.
Arguing as in \cite[Lemma 2]{horquint}
we get $2p_g-4\leq D^2\leq K_X^2$ and by
our assumptions we have equality everywhere.
Also, it is shown in (loc. cit.) that 
${\cal L}:=\omega_X\otimes\OO_D$ 
is an invertible
sheaf on $D$ with non-vanishing $h^0$ and $h^1$ for which Clifford's
inequality is an equality.
By Theorem \ref{clifford}, the curve $D$ is hyperelliptic 
in the sense of Definition \ref{hyperdef}.
We see ${\cal L}^{\otimes2}\iso\omega_D$ from the adjunction formula
on $X$
and so the map defined by $|{\cal L}|$ is not birational
by Lemma \ref{nonbirational}.
However, this contradicts the birationality of $\phi_1$.
Hence the canonical map cannot be birational.

Thus, $\deg\phi_1=2$ and the image of $\phi_1$ is an irreducible but
possibly singular surface
of degree at most $p_g-2$ in $\PP^{p_g-1}$.
However, $p_g-2$ is the lowest degree possible for a non-degenerate
surface in $\PP^{p_g-1}$ and so the image of $\phi_1$ has degree
equal to $p_g-2$.
We write $|K_X|=|M|+F$, where $F$ denotes the fixed part and $M$
denotes the movable part of the canonical system.
As $K_X$ is nef,
we have $2p_g-4=M^2\leq K_X^2$, and hence equality by our assumptions.
From $MF+F^2=K_X F\geq0$ and $MF\geq0$ together with
$K_X^2=M^2$ it is not difficult to deduce
$F^2=0$ and $K_XF=0$.
Since $K_X^2>0$ and $K_XF=F^2=0$, the Hodge index theorem implies
that $F$ is numerically trivial.
On the other hand, $F$ is an effective divisor and being numerically
trivial we see that $F$ is the zero divisor.\qed

\begin{Remark}
  A surface in $\PP^n$ that spans the ambient space has degree
  at least $n-1$.
  Surfaces of minimal degree have been classified by
  del~Pezzo (see \cite{eh} for a modern account) and consist of
  $\PP^2$, Hirzebruch surfaces, as well as cones over rational normal
  curves.
\end{Remark}

\section{Classification of Horikawa surfaces}
\label{classificationsection}

In this section we give a more detailed description of
Horikawa surfaces similar to 
\cite[Capitolo VIII.11]{enr} and \cite[Section 1]{horikawa}.
In order for this to work also in presence of inseparable maps
and wild ramification, we do not use the language of
branch divisors
but use the description of flat morphisms of
degree $2$ in terms of their associated line bundles.
\medskip

Let $\pi:X\to S$ be a flat double cover, i.e.,
a finite, flat and surjective 
morphism of degree $2$.
Via
\begin{equation}
\label{first}
0\,\to\,\OO_S\,\to\,\pi_\ast\OO_X\,\to\,{\cal L}^\vee\,\to\,0
\end{equation}
we obtain a sheaf ${\cal L}^\vee$ on $S$, which is invertible
by our assumptions on $\pi$.
In this case, we define $\cal L$ to be its dual.

\begin{Definition}
  We refer to $\cal L$ as the
  {\em line bundle associated with $\pi$}.
\end{Definition}

Flat double covers of smooth varieties are automatically
Gorenstein by \cite[Proposition 0.1.3]{cd}
and the dualising sheaf $\omega_X$ of $X$
is given by
\begin{equation}
\label{second}
\omega_X \,\iso\, \pi^\ast (\omega_S\otimes{\cal L})\,.
\end{equation}
Using the projection formula and (\ref{first}), we obtain
an extension
\begin{equation}
\label{pushforward}
0\,\to\,\omega_S\otimes{\cal L}\,\to\,
\pi_\ast\omega_X \,\to\,
\omega_S\,\to\,0\,.
\end{equation}
\begin{Remark}
  \label{canonicalfactorisation}
  If $p_g(S)=0$, then the canonical map of $X$
  factors as $\pi$ followed by the map from $S$ 
  associated with $\omega_S\otimes{\cal L}$. 
\end{Remark}
\medskip

For a natural number $d\geq0$, we denote by
$\FF_d$ the Hirzebruch surface 
$\PP_{\PP^1}(\OO_{\PP^1}\oplus\OO_{\PP^1}(d))$.
This $\PP^1$-bundle over $\PP^1$ has a section 
$\Delta_0$ with self-intersection number
$\Delta_0^2=-d$, which is unique if $d$ is positive.
We denote by $\Gamma$ the class of a fibre of this
$\PP^1$-bundle.

We now state and prove the structure result for
Horikawa surfaces.

\begin{Theorem}
  \label{classification}
  Let $X$ be a Horikawa surface and
  $S:=\phi_1(X)$ its canonical image in $\PP^{p_g-1}$.
  \begin{enumerate}
  \item If $S$ is smooth then $\phi_1$ factors as
     $$
     \phi_1\,:\, X\,\to\,\Xcan\,\stackrel{\pi}{\to}\,S\,,
     $$
     where $\Xcan$ denotes the canonical model of $X$ and
     where $\pi$ is a finite and flat morphism of degree $2$.
     Let $\cal L$ be the line bundle associated with $\pi$.
     Then we have the following possibilities:
     \begin{itemize}
     \item[-] $S\iso\PP^2$, $p_g=3$ and ${\cal L}\iso\OO_{\PP^2}(4)$,
     \item[-] $S\iso\PP^2$, $p_g=6$ and ${\cal L}\iso\OO_{\PP^2}(5)$,
     \item[-] $S\iso\FF_d$, 
     $0\leq d\leq p_g-4$, $p_g-d$ is even and 
     ${\cal L}\iso\OO_{\FF_d}(3\Delta_0+ \frac{1}{2}(p_g+2+3d)\Gamma)$.
     \end{itemize}
  \item If $S$ is not smooth then it is the cone over a rational
     normal curve of degree $p_g-2$ in $\PP^{p_g-1}$.
     Also $p_g\geq4$ and
     the minimal desingularisation $\nu:\tilde{S}\to S$ is 
     isomorphic to $\FF_{p_g-2}$.
     There exists a partial desingularisation $X'$
     of the canonical model $\Xcan$
     such that $\phi_1$ factors as
     $$
     \phi_1\,:\, X\,\to\,X'\,\stackrel{\pi}{\to}\,
     \FF_{p_g-2}\,\stackrel{\nu}{\to}\,S\,,
     $$
     where $\pi$ is a finite and flat morphism of degree $2$.
     If $\cal L$ denotes the line bundle associated with $\pi$ then
     ${\cal L}\iso\OO_{\FF_d}(3\Delta_0+ \frac{1}{2}(p_g+2+3d)\Gamma)$.
  \end{enumerate}
\end{Theorem}

\begin{Remark}
  We will see in Proposition \ref{restrictions} that 
  there are further restrictions on $p_g$ and $d$.
\end{Remark}

\prf
By construction, $\phi_1$ factors over the canonical
model $\Xcan$.
We know from Theorem \ref{firsttheorem} that $S$ is a surface
of degree $p_g-2$ in $\PP^{p_g-1}$, i.e., a surface of minimal
degree.
Such a surface is a smooth rational surface or the cone over
the rational normal curves of degree $p_g-2$ by
\cite{eh}.
Thus, if $S$ is not smooth its minimal desingularisation is 
the Hirzebruch surface $\FF_{p_g-2}$.

In any case, we denote by $\nu:\tilde{S}\to S$ the minimal
desingularisation of $S$.
Copying the proof of \cite[Lemma 1.5]{horikawa} we see
that $\phi_1$ factors over $\tilde{S}$.
Let $\psi$ be the induced morphism $X\to\tilde{S}$,
which is generically finite of degree $2$ by
Theorem \ref{firsttheorem}.
Let $X'$ be the Stein factorisation of $\psi$
and denote by $\pi$ the induced morphism $\pi:X'\to\tilde{S}$,
which is finite of degree $2$.
Moreover, $\pi$ is flat since $X'$ is Cohen--Macaulay
(being a normal surface)
and $\tilde{S}$ is regular.
Moreover, since flat double covers of smooth varieties
are Gorenstein, it follows that $X'$ is Gorenstein.

The morphism $X\to X'$ is birational.
Since $\phi_1$ factors over $\Xcan$, it is not difficult
to see that the canonical morphism $X\to\Xcan$ factors
over $X'$.
In particular, $X'$ has at worst rational singularities, which are
Gorenstein by what have already proved, i.e., $X'$ has
at worst Du~Val singularities.
As $\varphi:X'\to\Xcan$ is a birational morphism and $X'$
has only Du~Val singularities, $\Xcan$ is also the
canonical model of $X'$ and hence $\varphi$ partially 
resolves the singularities of $\Xcan$.

A smooth surface $S$ of minimal degree is either
$\PP^2$ embedded via $|H|$ or $|2H|$ into projective
space or a Hirzebruch surface $\FF_d$ embedded into
projective $n$-space via 
$|\Delta_0+\frac{1}{2}(n-1+d)\Gamma|$, where
$n-d-3$ is an even and non-negative integer
(see for example \cite[Lemma 1.2]{horikawa}).
Thus, if $X$ is a Horikawa surface and
$\phi_1(X)$ is a smooth surface $S$, the canonical
system of $X$ factors over $\omega_S\otimes\cal L$
on $S$, where $\omega_S\otimes\cal L$ is one of
the linear systems just described.
This yields the first list, cf. also 
\cite[Section 1]{horikawa}.

If $S$ is a singular surface of minimal degree,
then $\tilde{S}\iso\FF_{p_g-2}$ and the embedding
of $S$ is given by $|\Delta_0+d\Gamma|$ on $\tilde{S}$.
Proceeding as before, we obtain the description
of $\cal L$ in this case.
Again, we refer to \cite[Section 1]{horikawa} for
details.
\qed\medskip

Before proceeding we need (or recall)
a simple vanishing result.

\begin{Lemma}
 \label{technicalvanishing}
 We have $H^1(\PP^2,{\cal L})=0$ for every line bundle
 on $\PP^2$.

 On the Hirzebruch surface $\FF_d$ we have
 $H^1(\FF_d,\OO_{\FF_d}(a\Delta_0+b\Gamma))=0$ if
  \begin{enumerate}
  \item $a\geq0$ and $b\geq0$ or
  \item $a\leq-2$ and $b\leq-(d+2)$.
  \end{enumerate}
\end{Lemma}

\prf
We leave the assertion about line bundles on $\PP^2$
to the reader.

We consider the exact sequence
$$
0\,\to\,\OO_{\FF_d}((a-1)\Delta_0)\,\to\,\OO_{\FF_d}(a\Delta_0)
\,\to\,\OO_{\Delta_0}\,\to\,0\,.
$$
Taking cohomology and noting that the statement is clear
for $a=b=0$ we obtain the assertion for $a\geq0$ and $b=0$
inductively.
Using this result, an induction on $b$ shows
the vanishing for $a\geq0$ and $b\geq0$.
Applying Serre duality we obtain the
remaining vanishing result.
\qed\medskip

The following result is crucial to prove that 
Horikawa surfaces are simply connected as well as to show that 
there are further dependencies between 
$p_g$ and $d$ in Theorem \ref{classification}.

\begin{Lemma}
  \label{vanishing}
  Let $\cal L$ be as in Theorem \ref{classification}.
  Then $h^1({\cal L})=h^1({\cal L}^\vee)=0$.
\end{Lemma}

\prf
This follows immediately from inspecting the list of possible
${\cal L}$'s given in Theorem \ref{classification} and then
applying Lemma \ref{technicalvanishing}.
\qed\medskip

As an application of Theorem \ref{classification} we obtain

\begin{Proposition}
  \label{simplyconnected}
  A Horikawa surface fulfills $h^{01}(X):=h^1(\OO_X)=0$ and is 
  algebraically simply connected, i.e., has a trivial
  \'etale fundamental group.
  In particular, its Picard scheme is reduced.
\end{Proposition}

\prf
We use the notations of Theorem \ref{classification}.
If the canonical image $\phi_1(X)$ is a smooth surface, 
the long exact sequence of cohomology applied to
(\ref{first}) together with Lemma \ref{vanishing} 
yields $h^1(\OO_{\Xcan})=0$.
Since $\Xcan$ has at worst rational singularities,
the Grothendieck--Leray spectral sequence 
associated with the push-forward of the structure sheaf
yields $h^1(\OO_X)=0$.
The case where the canonical image $\phi_1(X)$ is
a singular surface is similar and left to the reader.

From $h^1(\OO_X)=0$ it follows that the Picard scheme 
of $X$ is reduced.

To prove that $X$ is algebraically simply connected we use
the idea of \cite[Theorem 14]{bom}:
Let $\hat{X}\to X$ be an \'etale cover of degree $m$.
Then we compute $\chi(\OO_{\hat{X}})=m\chi(\OO_X)$ and 
$K_{\hat{X}}^2=mK_X^2$.
Using Noether's inequality (Theorem \ref{noether}), 
we obtain
$$
  m\,(1+p_g(X)) \,=\, 1-h^1(\OO_{\hat{X}})+p_g(\hat{X})
  \,\leq\,1+p_g(\hat{X}) \,\leq\,
  \frac{1}{2}K_{\hat{X}}^2+3
  \,=\,\frac{m}{2}K_X^2+3\,.
$$
We assumed $X$ to be a Horikawa surface and so
this inequality holds for $m=1$ only.
Thus, every \'etale cover is trivial and hence $X$ is 
algebraically simply connected.
\qed

\begin{Remark}
  Over the complex numbers,
  even the topological fundamental group of a
  Horikawa surface is trivial
  \cite[Theorem 3.4]{horikawa}.
\end{Remark}
\medskip

Surfaces with $K^2=2$ and $p_g=3$, i.e., Horikawa
surfaces with $p_g=3$, have also been studied
in \v{S}afarevi\v{c}'s book \cite[Chapter 6.3]{saf}.
The emphasis there is on the fact that these complex surfaces
are the only surfaces besides those with $K^2=1$ and $p_g=2$
where $|3K_X|$ does not define a birational map.

\begin{Proposition}
  \label{non3birational}
  The $3$-canonical map of a minimal surface of general
  type with $K^2=2$ and $p_g=3$ is a morphism but
  not birational.
\end{Proposition}

\prf
From Theorem \ref{classification}, we see that the canonical
map of $X$ exhibits the canonical model $\Xcan$ as a flat double
cover $\pi:\Xcan\to\PP^2$ with associated line bundle 
${\cal L}\iso\OO_{\PP^2}(4)$.
Using (\ref{second}), we see that 
$\omega_{\Xcan}^{\otimes3}\iso\pi^\ast(\OO_{\PP^2}(3))$.
By the projection formula, the pushforward 
$\pi_\ast(\omega_{\Xcan}^{\otimes3})$ is an extension of
$\OO_{\PP^2}(-1)$ by $\OO_{\PP^2}(3)$.
Taking global sections, we see that
all global sections of $\omega_{\Xcan}^{\otimes3}$ are pull-backs
of global sections of $\OO_{\PP^2}(3)$.
Hence the $3$-canonical map of $X$ is a morphism and
factors over $\pi$.
In particular, it is not birational.
\qed

\section{The canonical double cover}

We now take a closer look at the canonical double
cover of a Horikawa surface.
In particular, we will see in 
Proposition \ref{restrictions} that
there are more restrictions on the line bundle $\cal L$
in Theorem \ref{classification}.
\medskip

Let $\pi:X\to S$ be a flat double cover
where $X$ is a normal and $S$ is a smooth variety.
Let $\cal L$ be the associated line bundle and consider
the short exact sequence (\ref{first}).
If the characteristic of the ground field is different from $2$ 
then the extension of function fields $k(X)/k(S)$ is Galois
with a cyclic Galois group of order $2$.
The Galois action extends to an action on $\pi_\ast\OO_X$
and decomposing into the $\pm1$-eigensheaves 
we obtain a splitting of (\ref{first}).

If the characteristic of the ground field is $p=2$ then
there is no reason for (\ref{first}) to split.
Of course, this sequence has to split if 
$ext^1({\cal L}^\vee,\OO_S)=h^1({\cal L})=0$.
By Lemma \ref{vanishing}, this condition is fulfilled for
Horikawa surfaces.

If $\pi$ is a flat double cover as above such that (\ref{first})
splits then there exist global sections $f\in H^0({\cal L})$
and $g\in H^0({\cal L}^{\otimes2})$ such that the cover
$\pi:X\to S$ is {\em globally} over $S$ given by
\begin{equation}
  \label{globaldoublecover}
  z^2\,+\,fz\,+\,g\,=\,0\,,
\end{equation}
where $z$ is a fibre coordinate on $\cal L$.
Moreover, if the characteristic $p$ is different from $2$ or if
$p=2$ and $\pi$ is purely inseparable, we may even assume
$f=0$.

\begin{Lemma}
  \label{reducedsections}
  Let $\pi:X\to S$ be a flat double cover
  where $X$ is normal and $S$ is smooth.
  In characteristic $p=2$ we assume moreover that the
  associated sequence (\ref{first}) splits.
  \begin{enumerate}
  \item If $p\neq2$ or $p=2$ and $\pi$ is purely inseparable
    then ${\cal L}^{\otimes2}$ has a reduced section.
  \item If $p=2$ and $\pi$ is separable then there exist 
    sections $f$ and $g$ of $\cal L$ and ${\cal L}^{\otimes2}$,
    respectively such that $\divisor(g)$ has no component of
    multiplicity $\geq2$ along a component of $\divisor(f)$.
  \end{enumerate}
\end{Lemma}

\prf
We may assume that $\pi$ is given by (\ref{globaldoublecover}).

In the first case we may assume $f=0$.
If $h^2$ divides $g$ then, by the Jacobian
criterion, $X$ is singular along the divisor
$\{\pi^\ast(h)=0\}$.
However, this contradicts the normality of $X$.

In the second case, if $h$ divides $f$ and 
$h^2$ divides $g$ then, again by the Jacobian criterion,
$X$ is singular along $\pi^\ast(h)=0$.
Again, this contradicts the normality of $X$.
\qed\medskip

This result imposes further restrictions on the line bundle
$\cal L$ of Theorem \ref{classification} and we obtain
the same result as in \cite[Section 1]{horikawa}, now valid
in arbitrary characteristic.

\begin{Proposition}
  \label{restrictions}
  Let $X$ be a Horikawa surface
  and $S$ be the image of its canonical map.
  \begin{enumerate}
  \item If $S$ is the smooth Hirzebruch surface $\FF_d$, then
  $p_g\geq2d-2$.
  \item If $S$ is singular, then $p_g\leq 6$.
  \end{enumerate}
\end{Proposition}

\prf
Let $\tilde{S}$ be the minimal desingularisation of $S$
and consider the finite and flat double cover
$\pi:X'\to\tilde{S}$ induced by $\phi_1$ as in 
Theorem \ref{classification}.
Let $\cal L$ be the line bundle associated with $\pi$
and $L$ be its class in $\NS(\tilde{S})$.
We know from Lemma \ref{vanishing} that $H^1({\cal L})$
vanishes and so the short exact sequence (\ref{first})
associated with $\pi$ splits.

We may assume that $d>0$ (else there is nothing to prove),
in which case the divisor $\Delta_0$ is unique.
We claim that $2L\cdot\Delta_0\geq-d$.
If not then every section of ${\cal L}^{\otimes2}$ vanishes
with multiplicity $\geq2$ along $\Delta_0$.
Hence every section of $\cal L$ vanishes
along $\Delta_0$.
This however, contradicts Lemma \ref{vanishing} in all
possible cases.

If $S=\FF_d$ then the inequality just shown together with
Theorem \ref{classification} yields $p_g\geq2d-2$.
In case $S$ is the cone over the rational normal curve of 
degree $d$ we obtain $p_g\leq6$.
\qed

\section{Inseparable canonical maps}
\label{inseparablesection}

Theorem \ref{classification} and 
Proposition \ref{restrictions} restrict the numerical
invariants of a Horikawa surface.
Conversely, we now establish the existence of Horikawa
surfaces for all remaining invariants.
If the characteristic is different from $2$, this can be
done along the classical arguments of Enriques and
Horikawa.
In characteristic $2$, we will first treat existence for
surfaces with inseparable canonical map.
In the next section we get the remaining case by
a deformation argument.
\medskip

A converse to Theorem \ref{classification} and 
Proposition \ref{restrictions}
in characteristic $p\neq2$ can be proved along the
lines of \cite[Section 1]{horikawa}.
We only sketch the argument here.

\begin{Theorem}[Horikawa]
  \label{horikawamainthm}
  Let $\tilde{S}$ be a minimal 
  rational surface in characteristic
  $p\neq2$ and $\cal L$ a line bundle on $\tilde{S}$.
  Assume moreover, that
  $\tilde{S}$ and $\cal L$ are as in Theorem \ref{classification}
  and that the additional inequalities of 
  Proposition \ref{restrictions} are fulfilled.
  Then there exists a Horikawa
  surface $X$ such that $\tilde{S}$ resolves the singularities of
  $\phi_1(X)$ and such that $\cal L$ is the line bundle associated
  with the flat double cover $\pi$ of Theorem \ref{classification}.
\end{Theorem}

\prf
Given a line bundle $\cal L$ on $\tilde{S}$ as above,
there exists a reduced section of 
${\cal L}^{\otimes 2}$ with at worst normal crossing singularities.
A double cover of $\tilde{S}$ branched along the divisor of this
section yields a surface with at worst Du~Val singularities, whose
desingularisation is a
Horikawa surface with the stated
properties, cf. \cite[Section 1]{horikawa}. 
\qed\medskip

Henceforth we shall assume that the characteristic of the ground
field is $p=2$.
To obtain examples with inseparable canonical maps,
we use quotients by foliations as described in 
\cite{lie}.
We note that this is a feature that does not occur for curves,
where the canonical map can never become inseparable.
Also, these surfaces of general type are automatically (inseparably)
unirational, which cannot happen in dimension $1$ 
by L\"uroth's theorem.

Although we only present examples we note that all Horikawa
surfaces with inseparable canonical map can be obtained this way.
In fact, the canonical morphism is inseparable if and only if
the morphism $\pi$ in Theorem \ref{classification} is purely 
inseparable.
Using the Frobenius morphism and Hirokado's result
\cite[Proposition 2.6]{hir1},  
that singularities of $p$-closed vector fields on surfaces
can be resolved by blow-ups
in characteristic $2$, 
we conclude that Horikawa surfaces
with inseparable canonical morphism are quotients of rational
surfaces by $p$-closed vector fields.
In particular, these surfaces are all (inseparably)
unirational.

\begin{Theorem}
  \label{smoothexistence}
  Given non-negative integers $g,d$ with
  $g\geq\max\{2d-2, d+4\}$ and $g-d$ even, 
  there exists
  a Horikawa surface with $p_g=g$ such that
  its canonical morphism maps inseparably onto
  $\FF_d$.  
\end{Theorem}

\prf
We use the coordinates of \cite{hir}
to express $\FF_d$ and its morphism onto $\PP^1$ as
$$
\left(\Proj k[X_1,Y_1]\times\Spec k[t_1]\right) 
\,\cup\,
\left(\Proj k[X_2,Y_2]\times\Spec k[t_2]\right) 
\,\to\,
\Proj k[T_1,T_2]\,.
$$
We set $x_i:=X_i/Y_i$ and $y_i:=Y_i/X_i$ for $i=1,2$.
Then
$$
\begin{array}{l@=ll@=l}
  x_1\, &\, x_2/t_2^d & \partial/\partial x_1\, &\, t_2^d\,\partial/\partial x_2\\
  y_i\, &\, 1/x_i & \partial/\partial x_i\, &\, -y_i^2\partial/\partial y_i\\
  t_1\, &\, 1/t_2 & \partial/\partial t_1\, &\, 
  -dx_2t_2\partial/\partial x_2 - t_2^2\partial/\partial t_2\,.
\end{array}
$$
The elements of $|\Gamma|$ are given by $\{t_i=const\}$
and the section $\Delta_0$ is defined by $\{y_1=0\}\cup\{y_2=0\}$,
whereas $\{x_1=0\}\cup\{x_2=0\}$ is linearly equivalent
to $\Delta_0+d\Gamma$.

We will first assume that $g$ and $d$ are even integers.
For an arbitrary integer $m\geq0$ we choose an integer
$\ell\geq0$ such that $0\leq 2\ell-2m+2\leq 5d$.
Then we choose pairwise distinct elements 
$a_1,...,a_\ell,b_1,...,b_m$ in the ground field $k$ 
and define the rational function  
$$
 \psi(t_1) \,:=\, \prod_{i=1}^\ell (t_1-a_i)^{-2}\,
 \prod_{j=1}^{m} (t_1-b_j)^{2}\,\in\,k(t_1)\,.
$$
Then we define the vector field
$$
\eta \,:=\, x_1^{-4}\frac{\partial}{\partial x_1} \,+\, 
\psi(t_1)\frac{\partial}{\partial t_1}
$$
which is additive in characteristic $2$, i.e., $\eta^{[2]}=0$.
Considered as a derivation of $k(x_1,t_1)$ over $k$ it extends to 
a rational vector field on $\FF_d$ with divisor
$$
(\eta) \,\lin\, -4\Delta_0 \,-\, (2m-2+4d)\Gamma\,.
$$
If $5d-2\ell+2m-2=0$ then $\eta$ has $\ell+m$ isolated singular points.
More precisely, the quotient $X':=\FF_d/\eta$ has $\ell$ singular points
of type $D_8$ and $m$ singular points of type $D_{12}$,
where we use \cite[Proposition 2.3]{lie} to determine the
singularities.
Using the results of \cite[Section 4]{lie} we see that $X'$
has an ample canonical sheaf if $g\geq d+4$, i.e., $X'$ is
the canonical model of a surface of general type.
We find $p_g=2m-2+2d\geq2d-2$ and that $X'$ is a Horikawa
surface.
Using Remark \ref{canonicalfactorisation}, we see that
the canonical map of $X'$ factors over $\FF_d^{(-1)}$, cf.
also \cite[Section 7]{lie}.
This construction yields all examples of Horikawa
surfaces with all possible values of $g$ and $d$ 
where both are even integers.

If $5d-2\ell+2m-2=4$ then the quotient $X':=\FF_d/\eta$ acquires
an elliptic singularity of type $(19)_0$ at the point lying below
$P:=\{x_2=t_2=0\}$ by \cite[Proposition 2.3]{lie}.
The other singular points of $X'$ are again Du~Val singularities
of type $D_8$ and $D_{12}$.
Desingularising $X'$ we obtain a Horikawa surface $X''$ with 
$p_g=2m-3+2d$.
By construction, the canonical image of $X''$ equals the
image of $\FF_d$ under the morphism that is defined
by imposing a simple base point at $P$ in the linear
system $|\Delta_0+\frac{1}{2}(g-2+d)\Gamma|$.
Blowing up $P$ resolves the indeterminacy and since $P$ does
not lie on $\Delta_0$, the induced morphism on the
blow-up of $\FF_d$ factors over $\FF_{d-1}$.
Hence the canonical system of $X''$ factors over $\FF_{d-1}$.
With this construction we obtain examples of Horikawa surfaces
with all possible values of $g$ and $d$ where both
are odd integers.
\qed

\begin{Theorem}
  \label{coneexistence}
  For every $4\leq g\leq 6$ there exists a Horikawa surface 
  with $p_g=g$ such that its canonical morphism maps 
  inseparably onto the cone over a rational normal curve
  of degree $p_g-2$ in $\PP^{p_g-1}$.
\end{Theorem}

\prf
To get an example with $p_g=4$ (resp. $p_g=6$) 
we take $d=2$ (resp. $d=4$), $m=1$ (resp. $m=0$) and 
$\ell=4$ (resp. $\ell=9$)
in the series of surfaces constructed in the proof of
Theorem \ref{smoothexistence}.

To get an example with $p_g=5$ we take $d=4$, $m=0$ 
and $\ell=7$.
The quotient $\FF_4/\eta$ has a singularity of type
$(19)_0$ and a desingularisation yields the desired surface.
\qed

\begin{Theorem}
  \label{existence}
  There exist Horikawa surfaces with $p_g=3$ and
  $K^2=2$ as well as $p_g=6$ and $K^2=8$ such that
  their canonical morphisms map inseparably onto
  $\PP^2$ and the Veronese surface in $\PP^5$,
  respectively.
\end{Theorem}

\begin{Remark}
  By a theorem of Bloch, there are no smooth and inseparable
  double covers of $\PP^2$, cf. \cite[Proposition 2.5]{ek1}.
  Hence the canonical model of a surface of
  Theorem \ref{existence} cannot be smooth. 
\end{Remark}

\prf
We use the Zariski surfaces constructed by Hirokado
\cite[Example 3.2]{hir}.
On $\FF_1$, the blow-up of $\PP^2$ in one point, the  
vector field $\Delta_1$ with $\ell=0$ and $n=3$ (resp. $n=4$) 
is multiplicative with isolated singularities of multiplicity $1$.
The quotient $X:=\FF_1/\Delta_1$ is a minimal surface of
general type with at worst $A_1$ singularities having the 
invariants we are looking for, cf.
\cite[Example 3.2]{hir} and \cite[Example 3.5]{hir}.

The canonical map of $X$ factors over $\FF_1^{(-1)}$ followed
by the map defined by the linear system
$|L+E|$ (resp. $|2L+2E|$), where $L$ is
the class of line pulled back from $\PP^2$ and $E$ is the
exceptional $(-1)$-curve on $\FF_1$.
Hence the canonical map of $X$ is a purely inseparable morphism 
onto $\PP^2$ (resp. onto the Veronese surface in $\PP^5$).
\qed

\begin{Remark}
 Alternatively, we can consider 
 $\delta:=(x^{2}+x^{-4})\frac{\partial}{\partial x}$, 
 which is a rational and additive vector field on $\PP^1$.
 The quotient $X':=(\PP^1\times\PP^1)/(\delta+\delta)$ has an elliptic
 singularity of type $(19)_0$ and nine Du~Val singularities of type
 $D_4$.
 Resolving these, we obtain another surface with $p_g=3$ and $K^2=2$
 such that the canonical map exhibits $\Xcan$ as a purely inseparable 
 double cover of $\PP^2$.
\end{Remark}

\section{Deformations to the separable case}
\label{deformationsection}

To obtain Horikawa surfaces with separable canonical map
in characteristic $2$, we
show that every Horikawa surface with inseparable canonical 
map can be deformed into such a surface
with a separable canonical map.
Morally speaking, the Horikawa surfaces with inseparable
canonical map should be boundary components of the moduli space
(if it exists in some sense)
of all Horikawa surfaces with fixed $p_g$ and fixed image
of the canonical map.
\medskip

First, we establish a result about deformations of
proper surface with at worst rational singularities.
If resolution of singularities in positive characteristic were
known to hold, we could probably argue along the lines of \cite{el} to
prove that every deformation of a rational singularity is again
rational.
For our purposes, it is enough to consider deformations over
a $1$-dimensional and normal base.
We prove the result in wider generality in order to apply it to
lifting problems later on as well.

\begin{Proposition}
  \label{deformationresult}
  Let $f:{\cal X}\to\Spec R$ be a flat and proper family of 
  surfaces where
  $R$ is a normal and $1$-dimensional Nagata ring. 
  Assume that the geometric fibre over some closed point 
  $t\in\Spec R$ is a 
  normal surface with at worst rational (resp. Du~Val) 
  singularities.
  Then there exists an open set $U\subseteq\Spec R$ containing 
  $t$ such that the geometric generic fibre as well as every 
  geometric fibre above points of $U$ is a 
  normal surface with at worst
  rational (resp. Du~Val) singularities.
\end{Proposition}

\prf
By \cite[Thm. 12.2.4]{ega4}, there exists an open and dense
subset $U$ of $C:=\Spec R$ over which the fibres are geometrically
normal.
Since $U$ and the fibres of $f$ above $U$
are normal then so is 
${\cal X}_U$, cf. \cite[Cor. 6.5.4]{ega4}.
Since any localisation of a normal ring remains normal, 
also the generic fibre ${\cal X}_K$ is normal.

The generic fibre ${\cal X}_K\to\Spec K$ is a 
surface over a field $K$ whose singularities
can be resolved by a
sequence of normalised blow-ups ${\cal Y}_K\to{\cal X}_K$
at closed points.
If ${\cal Y}_K$ is regular but not smooth over $K$ we pass
to a finite extension $L$ of $K$ where the singularities become
visible.
Base-changing to this field and resolving the singularities of ${\cal Y}_L$,
we eventually obtain a regular surface.
After possibly extending this field further, 
base-changing and resolving singularities we finally obtain a 
regular surface that is also smooth
over its ground field.
We may thus assume that there exists a finite
field extension $L$ of $K$ and a finite sequence of normalised
blow-ups ${\cal Y}_L\to{\cal X}_L$ such that ${\cal Y}_L$ is a
smooth surface over $L$.

Let $D$ be the normalisation of $C$ in $L$.
We choose a point $s\in D$ lying above $t\in C$.
Then the discrete valuation ring $\OO_{D,s}$ dominates
$\OO_{C,t}$. 
We denote their respective residue fields by $k$
and $\ell$.
Let ${\cal X}_D'$ be the normalisation of 
${\cal X}\otimes_{\OO_{C,t}}\OO_{D,t}$.
By the valuative criterion of properness the sequence
of normalised blow-ups $\nu_L:{\cal Y}_L\to{\cal X}_L$ 
extends to a sequence of normalised blow-ups 
$\nu_D':{\cal Y}_D'\to{\cal X}_D'$.
This induces a partial desingularisation of
the special fibre ${\cal X}_0\otimes_k \ell$.
Since this special fibre has at worst rational singularities,
the flat base-change theorem shows that 
$R^1(\nu_D')_\ast\OO_{{\cal Y}_D'}$ and
$R^1(\nu_L)_\ast\OO_{{\cal Y}_L}$ have to vanish.
In particular, ${\cal X}_L$ has
at worst rational singularities, which implies that
also the geometric generic fibre and
${\cal X}_K$ have at worst rational singularities.

There exists an open neighbourhood $V$ of $s\in D$ over which
${\cal Y}_D'\to{\cal X}_D'\to\Spec\OO_{D,s}$ 
spreads out flatly.
We have thus a sequence
of normalised blow-ups $\nu_V:{\cal Y}_V\to{\cal X}_V$ over $V$
which coincides with $\nu_L$ after tensorising with $L$.
Since ${\cal Y}_L$ is smooth over $L$, i.e., geometrically regular,
there exists a non-empty and open subset $W$ of $V$ over which
the fibres of ${\cal Y}_V\to V$ are geometrically regular,
cf. \cite[Thm. 12.2.4.]{ega4}.
Since $D$ is $1$-dimensional, the set $W\cup\{s\}$ is still open in
$D$ and we replace $V$ by it.
There exists a neighbourhood of $t$ such that the geometric
fibres of ${\cal X}_V\to V$ are normal (same argumentation as 
beginning of the proof) and we replace $V$ by it.
After possibly shrinking 
even further we may assume that 
$R^1(\nu_V)_\ast\OO_{{\cal Y}_V}$ is zero 
since this is true generically.
This open set $V$ now has the property that every geometric
fibre of ${\cal X}_V\to V$ is a normal surface with at worst 
rational singularities.

Since $R$ is a $1$-dimensional and normal Nagata ring
the map $D\to C$ is a finite morphism
between regular schemes and hence flat.
Thus, the image $W$ of $V$ in $C$ is open and contains
$t$.
We replace $U$ from above by
its intersection with $W$,
which yields an open set of $C$ containing $t$.
The geometric fibres of ${\cal X}\to C$ over $U$ 
are normal surfaces with at worst
rational singularities.

Now, suppose that the geometric fibre above 
$t$ has at worst Du~Val singularities.
These singularities are precisely the
rational Gorenstein singularities.
The special fibre is Gorenstein and since
the maximal ideal of $t\in C$ is generated by a regular
sequence, it follows that the generic fibre is Gorenstein
and there are points of 
$\cal X$ whose local rings are Gorenstein.
This argument also shows that the property of a fibre 
being Gorenstein is stable under generisation.
The set ${\cal S}$ of points of ${\cal X}$ such that
$\OO_{{\cal X},b}$ is not Gorenstein is
a proper closed subset of $\cal X$.
The image of $\cal S$ in $C$ is a constructible set and
hence so is its complement, which we have already seen
to be stable under generisation.
Hence that the set of points in $C$
whose fibres are Gorenstein schemes is open and 
contains $t$.
Intersecting $U$ with this open set we obtain an open
set of $C$ containing $t$ such that every 
geometric fibre above points of it
is a normal surface with at worst
rational Gorenstein, i.e., Du~Val, singularities.
\qed\medskip

The only point in the proof where we have used that
we deal with surfaces was when we used resolution of
singularities of the generic fibre.
For example, if the field of fractions of $R$ is of
characteristic zero we can apply Hironaka's
resolution of singularities and the proof works in arbitrary
dimensions.

Hence if we define a {\em weak rational singularity} of a variety $X$
in positive characteristic (naively) by requiring that
$R^if_\ast\OO_Y=0$ for every resolution $f:Y\to X$ and
every $i\geq1$ we obtain the

\begin{Corollary}
  Let $k$ be a perfect field of positive characteristic
  and $W(k)$ its associated Witt ring.
  Assume that $X$ is a normal variety over $k$ 
  with at worst weak rational (Gorenstein) singularities
  and that there exists a flat lifting 
  $f:{\cal X}\to\Spec W(k)$ of $X$.
  Then the generic fibre of $f$ 
  is a normal variety with at worst rational 
  (Gorenstein) singularities. \qed
\end{Corollary}
\medskip

We now apply the previous proposition to deform a given Horikawa
surface with inseparable canonical morphism into one with a separable 
canonical morphism.
This is achieved by deforming the canonical map while fixing
(the desingularisation of) the canonical image.

\begin{Theorem}
  \label{deformation}
  Let $X$ be a Horikawa surface in characteristic $2$ and assume
  that the canonical map $\phi_1:X\to S:=\phi_1(X)$ is inseparable.
  \begin{enumerate}
  \item If $S$ is smooth then the canonical model of $X$ can be deformed
     into the canonical model of a Horikawa surface with separable 
     canonical map and with the same canonical image.
  \item If $S$ is singular, then the surface $X'$ of 
     Theorem \ref{classification} can be deformed into 
     the $Y'$ of a Horikawa surface $Y$ with the same canonical image
     but with a separable canonical map.
  \end{enumerate}
\end{Theorem}

\prf
We do the case where $S$ is a smooth surface and leave the other
case to the reader.
Let $X$ be a Horikawa surface with inseparable canonical
map, $S$ its canonical image and
$\cal L$ be the line bundle associated with 
$\Xcan\to{S}$ as in Theorem \ref{classification}.
We know from Lemma \ref{vanishing} that
$\Xcan\iso\Spec{\cal A}$ with 
${\cal A}:=\OO_{S}\oplus{\cal L}$.
The algebra structure on $\cal A$ is given by
$z^2+t=0$ for some section $t$ of ${\cal L}^{\otimes 2}$.

We choose a non-zero global section $s$ of $\cal L$, which
exists by the description of $\cal L$ in 
Theorem \ref{classification}.
Then we define the 
$\OO_{S}$-algebra ${\cal A}_{\lambda}$ 
to be the $\OO_{{S}}$-module ${\cal A}$ with 
multiplication
$$
   z^2\,+\,{\lambda}sz\,+\,t\,=\,0\,.
$$
For $\lambda=0$ we obtain the original algebra $\cal A$.

We define $X_\lambda\,:=\,\Spec{\cal A}_\lambda$ and consider
these surfaces as a flat family of surfaces over the affine
line $\Aff^1_k$ with parameter $\lambda$.
Then, every surface in this family is a finite and 
flat double cover $\pi_\lambda:X_\lambda\to S$ 
with associated line bundle $\cal L$.
The fibre $X_0$ is the surface $X$ we started with.

There exists an open set $U\subseteq\Aff^1$ containing $\lambda=0$
such that $X_\lambda$ for $\lambda\in U$ is a normal surface with
at worst Du~Val singularities, cf. 
Proposition \ref{deformationresult}.
Using $\omega_{X_\lambda}\iso\pi_\lambda^\ast(\omega_S\otimes{\cal L})$
(cf. formula (\ref{second})) it is not difficult to see that the canonical
sheaves of these surfaces are ample, i.e., that all surfaces in
this family above $U$ are canonical models of Horikawa surfaces.
By Remark \ref{canonicalfactorisation}, their canonical morphisms 
factor as $\pi_\lambda$ followed by the morphism of $S$ associated
with the complete linear system $|{\cal L}|$.

As explained in the beginning of Section \ref{inseparablesection},
the zero set of $\lambda s$ determines the branch divisor of the
double cover $\pi_\lambda:X_\lambda\to\tilde{S}$.
If $\lambda\neq0$, then the morphism $\pi_\lambda$, i.e.,
the canonical morphism, is separable. 
\qed\medskip

We can summarise the results of Section \ref{inseparablesection} and 
Section \ref{deformationsection} as follows.

\begin{Theorem}
  Theorem \ref{horikawamainthm} also holds in characteristic $p=2$.
  Moreover, these surfaces exist with separable as well as inseparable
  canonical maps.\qed
\end{Theorem}

\section{Lifting to characteristic zero}

We now prove that Horikawa surfaces lift projectively over
$W(k)$.
By Theorem \ref{classification}, giving a Horikawa surface
is the same thing as giving a rational surface, a line bundle
on it and two sections.
It is not difficult to see that all this data lifts to
characteristic zero from which we obtain a lift of the canonical model of
any given Horikawa surface.
To achieve a lifting of the desingularisation we use Artin's result
on simultaneous resolution of singularities, which provides us with
an algebraic space over a possibly ramified extension of $W(k)$
achieving the lift.

Let $k$ be an algebraically closed field of 
characteristic $p>0$ and $W(k)$ its associated 
Witt ring.
The following results is probably folklore but
worthwhile being stated explicitly.

\begin{Lemma}
  \label{lifting}
  Let $S$ be a smooth rational surface over $k$ and
  $\cal L$ a line bundle on $S$.
  Then $S$ lifts to a surface over $W(k)$ and 
  $\cal L$ lifts to a unique line bundle $\widetilde{\cal L}$
  on this lift.
  
  In case $h^1({\cal L})=0$ we have furthermore 
  $h^i({\cal L})=h^i(\widetilde{{\cal L}})$ for all $i$.
  In particular, sections of $\cal L$ lift as well in this case.
\end{Lemma}

\prf
The first assertions follow from Grothendieck's existence
theorem and elementary deformation theory.
That smooth rational surfaces lift projectively over $W(k)$
is for example explained in \cite[Section 8.5.26]{ill2}.
From \cite[Corollary 8.5.6]{ill2} and 
$h^1(S,\OO_S)=h^2(S,\OO_S)=0$ we get the unique 
lifting of line bundles.

Let $\widetilde{{\cal L}}$ be the unique lift of $\cal L$ on $S$ 
and assume that $h^1({\cal L})=0$.
From the upper semicontinuity theorem we get $h^1(\tilde{{\cal L}})=0$.
Since $\chi({\cal L})=\chi(\widetilde{{\cal L}})$ and $h^1$ of both
line bundles is zero,
we conclude $h^i({\cal L})=h^i(\widetilde{{\cal L}})$ for $i=0,2$,
again using upper semicontinuity.
\qed\medskip

We now come to the main result of this section.

\begin{Theorem}
  Let $X$ be a Horikawa surface over a field $k$
  of positive characteristic.
  Then the canonical model of $X$ lifts over
  the Witt ring $W(k)$.
  Also, $X$ can be lifted in the category of algebraic spaces,
  i.e., there exists an algebraic space, flat over a possibly ramified 
  extension of $W(k)$, with special fibre $X$.
\end{Theorem}

\prf
We first do the case 
where the canonical image $S$ of $X$ is a smooth 
surface.

Let $\pi:\Xcan\to S$ the canonical flat double cover
with associated line bundle $\cal L$ as in 
Theorem \ref{classification}.
As an $\OO_S$-algebra $\Xcan$ is isomorphic
to ${\cal A}:=\OO_S\oplus{\cal L}$ by
Lemma \ref{vanishing}.
The map $\pi$ is globally given 
in the form (\ref{globaldoublecover})
for sections $s,t$ of
$\cal L$, ${\cal L}^{\otimes2}$, respectively.

By Lemma \ref{lifting}, we can lift $S$,
$\cal L$ and ${\cal L}^{\otimes2}$ over $W(k)$.
Let ${\cal S}\to W(k)$ be the lift of $S$ 
and ${\cal S}_K$ the generic fibre.
Using Lemma \ref{vanishing} and Lemma \ref{lifting}
again, we also lift the sections $s$ and $t$.
Out of this data we construct a flat double
${\cal X}'\to{\cal S}$ with special fibre 
$\Xcan\to S$.
By Proposition \ref{deformationresult} the generic
fibre ${\cal X}_K'$ 
of ${\cal X}'\to W(k)$ is a normal
surface with at worst Du~Val singularities.
Using the explicit description it follows
that ${\cal X}_K$ is a 
possibly singular Horikawa surface 
with at worst Du~Val singularities
and with canonical image ${\cal S}_K$.
Hence, we have lifted the canonical model of
$X$ over $W(k)$.

We have to resolve the singularities in the
special fibre.
However,  to achieve a simultaneous resolution
of singularities, we may have to base-change to an
algebraic space of finite type over $W(k)$,
cf. \cite[Theorem 1]{artin}.
That this resolution is in fact possible after 
a finite extension of $W(k)$ 
follows from the fact that $W(k)$ is Henselian, 
cf. \cite[Theorem 2]{artin}.

In case the canonical image is singular we proceed
as before to obtain a lift of $X'$ (in the 
notation of Theorem \ref{classification}) over $W(k)$ 
and to get a lift of $X$ over a possibly ramified 
extension of $W(k)$.
The pluri-canonical ring associated with the lift of
$X'$ specialises to the canonical model of $X$ and
generalises to a normal surface with at worst Du~Val
singularities.
This achieves a lifting of the canonical model of
$X$ also in the case of a singular canonical image.
\qed


\begin{thebibliography}{EGA IV}
  \bibitem[Ar]{artin} M.~Artin, {\it Algebraic Construction of Brieskorn's 
    Resolutions}, J. of Algebra 29, 330-348 (1974).
  \bibitem[BHPV]{bhpv} W.~P.~Barth, K.~Hulek, C.~Peters, A. van de Ven,
     {\it Compact Complex Surfaces}, 2nd edition, Erg. d. Math., 3. Folge
     Volume 4, Springer (2004). 
  \bibitem[Bo]{bom} E.~Bombieri, {\it Canonical models of surfaces of general
     type}, Inst. Hautes \'Etudes Sci. Publ. Math. 42, 171-220 (1973).
  \bibitem[CD]{cd} F.R.~Cossec, I.V.~Dolgachev, {\it Enriques Surfaces I},
     Prog. in Math. 76, Birkh\"auser (1989).
  \bibitem[EH]{eh} D.~Eisenbud, J.~Harris, {\it On Varieties of Minimal Degree
     (A Centennial Account)}, Algebraic Geometry, Bowdoin 1985, Proc. Symp. Pure
     Math. 46, Part 1, 3-13 (1987).
  \bibitem[EKS]{eks} D.~Eisenbud, J.~Koh, M.~Stillman, {\it Determinantal 
     equations for curves of high degree},
     Am. J. Math. 110, 513-539 (1988).
  \bibitem[E]{ek1} T.~Ekedahl, {\it Foliations and Inseparable Morphisms}, 
     Proc. Symp. Pure Math. 46, No. 2, 139-149 (1987).
  \bibitem[E2]{ek2} T.~Ekedahl, {\it Canonical models of surfaces of general type in 
     positive characteristic},  Inst. Hautes \'Etudes Sci. Publ. Math.  No. 67, 
     97-144 (1988).
  \bibitem[El]{el} R.~Elkik, {\it Singularit\'es rationelles et d\'eformations},
     Invent. math. 47, 139-147 (1978).
  \bibitem[En]{enr} F.~Enriques, {\it Le superficie algebriche},
     Nicola Zanichelli (1949).
  \bibitem[EGA IV]{ega4} A.~Grothendieck, {\it \'El\'ements de g\'eom\'etrie 
     alg\'ebrique IV: \'Etude locale des sch\'emas et des morphismes de
     sch\'emas}, Publ. Math. IHES 20 (1964), 24 (1965), 28 (1966), 32 (1967).
  \bibitem[Hart]{hart} R.~Hartshorne, {\it Algebraic Geometry}, Springer (1977).
  \bibitem[Hir1]{hir1} M.~Hirokado, {\it Singularities of multiplicative $p$-closed
     vector fields and global $1$-forms on Zariski surfaces}, J.~Math.~Kyoto Univ.,
     39-3, 455-468 (1999)
 \bibitem[Hir2]{hir} M.~Hirokado, {\it Zariski surfaces as quotients of Hirzebruch
     surfaces by 1-foliations}, Yokohama Math. J. 47, 103-120 (2000).
  \bibitem[Hor1]{horquint} E.~Horikawa, {\it On Deformations of Quintic Surfaces},
     Invent.~math. 31, 43-85 (1975).
  \bibitem[Hor2]{horikawa} E.~Horikawa, {\it Algebraic surfaces of general
     type with small $c_1^2$, I}, Ann. Math. 104, 357-387 (1976).
  \bibitem[Ill]{ill2} L.~Illusie, {\it Grothendieck's existence theorem in formal
    geometry with a letter of Jean-Pierre Serre}, 
    AMS Math. Surveys Monogr. 123, Fundamental Algebraic Geometry, 
    179-233 (2005). 
  \bibitem[Jou]{jou} J.-P.~Jouanolou, {\it Th\'eor\`emes de Bertini et
     Applications}, Prog. in Math. 42, Birkh\"auser (1983).
  \bibitem[Lie]{lie} C.~Liedtke, {\it Uniruled Surfaces of General Type},
     arXiv:math.AG/0608604 (2006), to appear in Math.Z.
  \bibitem[Noe]{no} M.~Noether, {\it Zur Theorie der eindeutigen Entsprechungen
      algebraischer Gebilde}, Math. Ann. 8, 495-533 (1875).
   \bibitem[S]{saf} I.~R.~\v{S}afarevi\v{c} et al., {\it Algebraic Surfaces},
      Proc. Steklov Inst. 75 (1965).
\end{thebibliography}
\end{document}